\documentclass[11pt]{amsart}
\usepackage{amsmath,amsthm,amsfonts,amssymb,latexsym,epsfig}

\newtheorem{theorem}{Theorem}[section]

\numberwithin{equation}{section}

\theoremstyle{definition}

\theoremstyle{remark}

\begin{document}
\title{On A Mixed Arithmetic-Geometric Mean Inequality}
\author{Peng Gao}
\address{Department of Mathematics, School of Mathematics and System Sciences, Beijing University of Aeronautics and Astronautics, P. R. China}
\email{penggao@buaa.edu.cn}

\begin{abstract}We extend a result of Holland on a mixed arithmetic-geometric mean inequality.
\end{abstract}

\maketitle
\section{Introduction}
\label{sec1}
   Let $M_{n,r}({\bf q},{\bf x})$ be the generalized weighted power means:
   $M_{n,r}({\bf q}, {\bf x})=(\sum_{i=1}^{n}q_ix_i^r)^{\frac {1}{r}}$,
   where ${\bf q}=(q_1,q_2,\cdots,
   q_n)$, ${\bf x}=(x_1,x_2,\cdots,
   x_n)$, $q_i >0, 1 \leq i \leq n$ with $\sum_{i=1}^nq_i=1$. Here
   $M_{n,0}({\bf q}, {\bf x})$ denotes the limit of $M_{n,r}({\bf q}, {\bf x})$ as
   $r\rightarrow 0^{+}$. Unless specified, we always assume $x_i>0, 1 \leq i \leq n$.
    When there is no risk of confusion,
    we shall write $M_{n,r}$ for $M_{n,r}({\bf q}, {\bf x})$ and
    we also denote $A_n, G_n$ for the arithmetic mean $M_{n,1}$, geometric mean $M_{n,0}$, respectively.

   For fixed ${\bf x}=(x_1,\cdots, x_n),{\bf w}=(w_1, \cdots, w_n)$ with $w_1>0, w_i \geq 0$, we define ${\bf x}_i=(x_1, \cdots, x_i)$, ${\bf w}_i=(w_1,
\cdots, w_i), W_i=\sum^i_{j=1}w_j$, $M_{i,
   r}=M_{i,r}({\bf w}_i/W_i, {\bf x}_i), {\bf M}_{i,r}=(M_{1,r}, \cdots,
M_{i,r})$. The following result on mixed mean inequalities is due to Nanjundiah \cite{N} (see also \cite{B}):
\begin{theorem}
\label{thm1.0}
  Let $r >s$ and $n \geq 2$. If for $
   2 \leq k \leq n-1$, $W_nw_k-W_kw_n>0$. Then
\begin{equation*}
   M_{n,s}({\bf M}_{n,r}) \geq M_{n,r}(
   {\bf M}_{n,s}),
\end{equation*}
   with equality holding if and only if $x_1 = \cdots =x_n$.
\end{theorem}

  It is easy to see that the case $r=1, s=0$ of Theorem \ref{thm1.0} follows from the following Popoviciu-type inequalities established in \cite{Ked1} (see also \cite[Theorem 9]{B}):
\begin{theorem}
\label{thm2}
  Let $n \geq 2$. If for $2 \leq k \leq n-1$, $W_nw_k-W_kw_n>0$, then
\begin{equation*}
   W_{n-1}\Big(\ln M_{n-1,0}({\bf M}_{n-1,1})-\ln M_{n-1,1}({\bf M}_{n-1,0}) \Big) \leq W_n\Big (\ln M_{n,0}({\bf M}_{n,1})-\ln M_{n,1}({\bf M}_{n,0}) \Big)
\end{equation*}
   with equality holding if and only if $x_n=M_{n-1,0}=M_{n-1,1}({\bf M}_{n-1,0})$.
\end{theorem}
   In \cite{T&T}, the following Rado-type inequalities were established:
\begin{theorem}
\label{thm1'}
  Let $s<1 $ and $n \geq 2$. If for $
   2 \leq k \leq n-1$, $W_nw_k-W_kw_n>0$, then
\begin{equation*}
   W_{n-1}\Big (M_{n-1,s}({\bf M}_{n-1,1})-M_{n-1,1}({\bf M}_{n-1,s}) \Big ) \leq W_n \Big (M_{n,s}({\bf M}_{n,1})-M_{n,1}({\bf M}_{n,s}) \Big )
\end{equation*}
   with equality holding if and only if $x_1 = \cdots =x_n$ and the above
   inequality reverses when $s>1$.
\end{theorem}

   The above theorem is readily seen to imply Theorem \ref{thm1.0}. In \cite{H}, Holland further improved the condition in Theorem \ref{thm1'} for the case $s=0$ by proving the following:
\begin{theorem}
\label{thm3}
  Let $n \geq 2$. If $W^2_{n-1} \geq w_{n}\sum^{n-2}_{i=1}W_i$ with the empty sum being $0$, then
\begin{equation}
\label{1.3}
   W_{n-1}\Big (M_{n-1,0}({\bf M}_{n-1,1})-M_{n-1,1}({\bf M}_{n-1,0}) \Big ) \leq W_n \Big (M_{n,0}({\bf M}_{n,1})-M_{n,1}({\bf M}_{n,0}) \Big )
\end{equation}
   with equality holding if and only if $x_1 = \cdots =x_n$.
\end{theorem}

  It is our goal in this paper to extend the above result of Holland by considering the validity of inequality \eqref{1.3} for the case $W^2_{n-1} < w_{n}\sum^{n-2}_{i=1}W_i$. Note that this only happens when $n \geq 3$. In the next section, we apply the approach in \cite{G} to prove the following
\begin{theorem}
\label{thm4}
  Let $n \geq 3$.  Inequality \eqref{1.3} holds when the following conditions are satisfied:
\begin{align}
\label{2.10}
 & 0 < \frac {w_{n}\sum^{n-2}_{i=1}W_i}{W^2_{n-1}}-1  \leq \frac {w_1}{w_n}, \quad \frac {W_{n-1}}{W_n}\prod^{n-2}_{i=1}\left ( \frac {W_{i+1}}{W_i} \right )^{\frac {W_iw_n}{W^2_{n-1}}} \leq 1, \\
 & \left (\frac {W_{n-1}w_{n}}{W_nw_1}\left (\sum^{n-2}_{i=1} \frac {W_iw_n}{W^2_{n-1}}-1 \right )+ \frac {w_{n}}{W_n} \right )\prod^{n-2}_{i=1}\Big(\frac {W_{i+1}}{w_{i+1}}\Big )^{\frac {w_{i+1}}{W_{n-1}}} \leq 1. \nonumber
\end{align}
\end{theorem}

  In the above theorem, we do not give the condition for the equality in \eqref{1.3} to hold. In Section \ref{sec 3}, we show that there do exist sequences $\{ w_i \}^n_{i=1}$  that satisfy the conditions of Theorem \ref{thm4}. 

\section{Proof of Theorem \ref{thm4}}
\label{sec 2} \setcounter{equation}{0}
    We may assume that $x_i>0, w_i > 0, 1 \leq i \leq n$ and the case $x_i=0$ or $w_i=0$ for some $i$ follows by continuity. We recast \eqref{1.3} as
\begin{equation}
\label{2.3}
   G_{n}({\bf A}_{n})-\frac {W_{n-1}}{W_n}G_{n-1}({\bf A}_{n-1})-\frac {w_n}{W_n}G_n \geq 0.
\end{equation}
  Note that
\begin{align}
\label{1.6}
  G_{n}({\bf A}_{n}) &=& \Big(G_{n-1}({\bf A}_{n-1}) \Big )^{W_{n-1}/W_n}A^{w_{n}/W_n}_n, \quad G_{n-1}({\bf A}_{n-1}) &=& A_n\prod^{n-1}_{i=1}\Big(\frac {A_i}{A_{i+1}}\Big )^{W_i/W_{n-1}}.
\end{align}
  Dividing $G_{n}({\bf A}_{n})$ on both sides of \eqref{2.3} and using \eqref{1.6}, we can recast \eqref{2.3} as:
\begin{equation}
\label{1.7}
   \frac {W_{n-1}}{W_n}\prod^{n-1}_{i=1}\Big(\frac {A_i}{A_{i+1}}\Big )^{W_iw_n/(W_{n-1}W_n)}+\frac {w_n}{W_n}\prod^{n}_{i=1}\Big(\frac {x_i}{A_{i}}\Big )^{w_i/W_n} \leq 1.
\end{equation}
   We express $x_i=(W_iA_i-W_{i-1}A_{i-1})/w_i, 1\leq i \leq n$ with $W_0=A_0=0$ to recast \eqref{1.7} as
\begin{equation*}
   \frac {W_{n-1}}{W_n}\prod^{n-1}_{i=1}\Big(\frac {A_i}{A_{i+1}}\Big )^{W_iw_n/(W_{n-1}W_n)}+\frac {w_n}{W_n}\prod^{n}_{i=1}\Big(\frac {W_iA_i-W_{i-1}A_{i-1}}{w_iA_{i}}\Big )^{w_i/W_n} \leq 1.
\end{equation*}
  We set $y_i=A_{i}/A_{i+1}$, $1 \leq i \leq 2$ to further recast the above inequality as
\begin{equation}
\label{1.8}
  \frac {W_{n-1}}{W_n}\prod^{n-1}_{i=1}y^{W_iw_n/(W_{n-1}W_n)}_i+\frac {w_n}{W_n}\prod^{n-1}_{i=1}\Big(\frac {W_{i+1}}{w_{i+1}}-\frac {W_{i}}{w_{i+1}}y_i\Big )^{w_{i+1}/W_n} \leq 1.
\end{equation}
 
   We now regard the right-hand side expression above as a function of $y_{n-1}$ only and define
\begin{align*}
   f(y_{n-1})=\frac {W_{n-1}}{W_n}c \cdot y^{w_n/W_n}_{n-1}+\frac {w_n}{W_n}c'\cdot \Big(\frac {W_{n}}{w_{n}}-\frac {W_{n-1}}{w_{n}}y_{n-1}\Big )^{w_n/W_n},
\end{align*}
  where
\begin{align*}
  c=\prod^{n-2}_{i=1}y^{W_iw_n/(W_{n-1}W_n)}_i, \quad c'=\prod^{n-2}_{i=1}\Big(\frac {W_{i+1}}{w_{i+1}}-\frac {W_{i}}{w_{i+1}}y_i\Big )^{w_{i+1}/W_n}.
\end{align*}
   On setting $f'(y_{n-1})=0$, we find that
\begin{align*}
   y_{n-1}=\left (\frac {W_{n-1}}{W_{n}}+\frac {w_{n}}{W_{n}}\left (\frac {c'}{c} \right )^{W_n/W_{n-1}} \right )^{-1}.
\end{align*}
   It is easy to see that $f(y_{n-1})$ is maximized at the above value with its maximal value being
\begin{align*}
  \left (\frac {W_{n-1}}{W_{n}}c^{W_n/W_{n-1}}+\frac {w_{n}}{W_{n}}(c')^{W_n/W_{n-1}} \right)^{W_{n-1}/W_{n}}.
\end{align*}
   Thus, in order for inequality \eqref{1.8} to hold, it suffices to have
\begin{align*}
  \frac {W_{n-1}}{W_{n}}c^{W_n/W_{n-1}}+\frac {w_{n}}{W_{n}}(c')^{W_n/W_{n-1}} \leq  1. 
\end{align*}
   Explicitly, the above inequality is
\begin{align*}
  g(y_1, y_2, \ldots, y_{n-2}):=\frac {W_{n-1}}{W_{n}}\prod^{n-2}_{i=1}y^{W_iw_n/W^2_{n-1}}_i+\frac {w_{n}}{W_{n}}\prod^{n-2}_{i=1}\Big(\frac {W_{i+1}}{w_{i+1}}-\frac {W_{i}}{w_{i+1}}y_i\Big )^{w_{i+1}/W_{n-1}} \leq  1. 
\end{align*}

 Let $(a_1, a_2, \ldots, a_{n-2}) \in [0, W_2/W_1]\times [0, W_3/W_2]\times \ldots \times [0, W_{n-1}/W_{n-2}]$ be the point in which the absolute maximum of $g$ is reached. If one of the $a_i$ equals $0$ or $W_{i+1}/W_i$, then it is easy to see that we have
\begin{align}
\label{25}
   g(a_1, a_2, \ldots, a_{n-2}) \leq \max \left ( \frac {W_{n-1}}{W_n}\prod^{n-2}_{i=1}\left ( \frac {W_{i+1}}{W_i} \right )^{W_iw_n/W^2_{n-1}}, \quad \frac {w_n}{W_n}\prod^{n-2}_{i=1}\left (\frac {W_{i+1}}{w_{i+1}} \right )^{w_{i+1}/W_{n-1}} \right ).
\end{align}
   If the point  $(a_1, a_2, \ldots, a_{n-2})$ is an interior point, then we have
\begin{align*}
  \nabla g(a_1, a_2, \ldots, a_{n-2})=0.
\end{align*}
   It follows that for every $1 \leq i \leq n-2$, we have
\begin{align}
\label{2.05}
 \frac {\prod^{n-2}_{i=1}a^{W_iw_n/W^2_{n-1}}_i}{\prod^{n-2}_{i=1}\Big(\frac {W_{i+1}}{w_{i+1}}-\frac {W_{i}}{w_{i+1}}a_i\Big )^{w_{i+1}/W_{n-1}}}=\frac {a_i}{\frac {W_{i+1}}{w_{i+1}}-\frac {W_{i}}{w_{i+1}}a_i}:=\frac 1d, 
\end{align}
  where $d>0$ is a constant (depending on the $w_i$). In terms of $d$, we have
\begin{align*}
   a_i=\frac {W_{i+1}}{dw_{i+1}+W_i}.
\end{align*}
    We use this to recast the first equation in \eqref{2.05} as
\begin{align*}
  \prod^{n-2}_{i=1}\left(\frac {W_{i+1}}{dw_{i+1}+W_i} \right )^{W_iw_n/W^2_{n-1}}=\frac 1d\prod^{n-2}_{i=1}\Big(\frac {W_{i+1}}{w_{i+1}}\cdot \frac {dw_{i+1}}{dw_{i+1}+W_i}\Big )^{w_{i+1}/W_{n-1}}.
\end{align*}
   We recast the above equality as
\begin{align*}
 &\ln \left ( \prod^{n-2}_{i=1}\left ( W_{i+1} \right )^{W_iw_n/W^2_{n-1}-w_{i+1}/W_{n-1}} \right )  \\
= & \sum^{n-2}_{i=1}\left ( \frac {W_iw_n}{W^2_{n-1}}-\frac {w_{i+1}}{W_{n-1}} \right )\ln \left ( dw_{i+1}+W_i \right )-\frac {w_1}{W_{n-1}}\ln d:=h(d).
\end{align*}
   Note that the above equality holds when $d=1$ and we have
\begin{align*}
   h'(d)=\sum^{n-2}_{i=1}\left ( \frac {W_iw_n}{W^2_{n-1}}-\frac {w_{i+1}}{W_{n-1}} \right )\frac 1{d+W_i/w_{i+1}} -\frac {w_1}{W_{n-1}}\frac 1d.
\end{align*}
   As 
\begin{align*}
   \frac {W_iw_n}{W^2_{n-1}}-\frac {w_{i+1}}{W_{n-1}}  \geq 0  \Leftrightarrow \frac {W_i}{w_{i+1}} \geq \frac {W_{n-1}}{w_n},
\end{align*}
   it follows that
\begin{align*}
   h'(d) & \leq \sum^{n-2}_{i=1}\left ( \frac {W_iw_n}{W^2_{n-1}}-\frac {w_{i+1}}{W_{n-1}} \right )\frac 1{d+W_{n-1}/w_{n}} -\frac {w_1}{W_{n-1}}\frac 1d \\
   & \leq  \frac {\left (\sum^{n-2}_{i=1} \frac {W_iw_n}{W^2_{n-1}}-1 \right )d-\frac {w_1}{w_n}}{d(d+W_{n-1}/w_{n})}.
\end{align*}
   Therefore, when
\begin{align*}
  \sum^{n-2}_{i=1} \frac {W_iw_n}{W^2_{n-1}}-1  \leq 0,
\end{align*}
   the function $h(d)$ is a decreasing function of $d$ so that $d=1$ is the only value that satisfies \eqref{2.05} and we have $a_i=1$ correspondingly with $g(1, 1, \ldots, 1)=1$ and this allows us to recover Theorem \ref{thm3}, by combining the observation that the right-hand side expression of \eqref{25} is an increasing function of $w_n$ for fixed $w_i, 1 \leq i \leq n-1$ with the discussion in the next section.

  Suppose now 
\begin{align}
\label{2.08}
  \sum^{n-2}_{i=1} \frac {W_iw_n}{W^2_{n-1}}-1  > 0,
\end{align}
  then the function $h(d)$ is a decreasing function of $d$ for
\begin{align*}
  d \leq \frac {\frac {w_1}{w_n}}{\sum^{n-2}_{i=1} \frac {W_iw_n}{W^2_{n-1}}-1}:=d_0.
\end{align*}
   If follows that if $d_0 \geq 1$, then  $d=1$ is the only value $\leq d_0$ that satisfies \eqref{2.05} and we have $a_i=1$ correspondingly with $g(1, 1, \ldots, 1)=1$. We further note that for any $d \geq d_0$ satisfying \eqref{2.05}, the value of $g$ at the corresponding $a_i$ satisfies 
\begin{align*}
   g(a_1, a_2, \ldots, a_{n-2})&= \left (\frac {W_{n-1}}{dW_n}+ \frac {w_{n}}{W_n} \right )\prod^{n-2}_{i=1}\Big(\frac {W_{i+1}}{w_{i+1}}\cdot \frac {dw_{i+1}}{dw_{i+1}+W_i}\Big )^{w_{i+1}/W_{n-1}} \\
& \leq \left (\frac {W_{n-1}}{dW_n}+ \frac {w_{n}}{W_n} \right )\prod^{n-2}_{i=1}\Big(\frac {W_{i+1}}{w_{i+1}}\Big )^{w_{i+1}/W_{n-1}} \\
& \leq \left (\frac {W_{n-1}}{d_0W_n}+ \frac {w_{n}}{W_n} \right )\prod^{n-2}_{i=1}\Big(\frac {W_{i+1}}{w_{i+1}}\Big )^{w_{i+1}/W_{n-1}}.
\end{align*}

   Combining this with \eqref{25}, we see that inequality \eqref{1.8} holds when the conditions in \eqref{2.10} are satisfied and this completes the proof of Theorem \ref{thm4}.

\section{A Further Discussion}
\label{sec 3} \setcounter{equation}{0}
  We show in this section that there does exist sequences $\{ w_i \}^n_{i=1}$ satisfying the conditions of Theorem \ref{thm4}. To see this, we note that the left-hand side expression of \eqref{2.08} vanishes when 
\begin{align}
\label{2.11}
  w_n=\frac {W^2_{n-1}}{\sum^{n-2}_{i=1}W_i}.
\end{align}
   
  It follows by continuity that such sequences  $\{ w_i \}^n_{i=1}$ satisfying the conditions of Theorem \ref{thm4} exist as long as the positive sequence  $\{ w_i \}^n_{i=1}$ with $w_i, 1 \leq i \leq n-1$ being arbitrary and  $w_n$ defined by \eqref{2.11} satisfies the last two inequalities of \eqref{2.10} with strict inequalities there. It is readily checked that these inequalities become
\begin{align}
\label{2.12}
  & \left ( \frac {\sum^{n-2}_{i=1}W_i}{\sum^{n-1}_{i=1}W_i} \right )^{\sum^{n-2}_{i=1}W_i}W^{W_{n-2}}_{n-1}<\prod^{n-2}_{i=1}W^{w_i}_i, \\ 
\label{2.13}
  & \frac {W_{n-1}}{\sum^{n-1}_{i=1}W_i} \prod^{n-2}_{i=1}\Big(\frac {W_{i+1}}{w_{i+1}}\Big )^{w_{i+1}/W_{n-1}}<1.
\end{align}

    Now, it is easy to see that inequality \eqref{2.12} holds when $n=3$. It follows that it holds for all $n \geq 3$ by induction as long as we have
\begin{align}
\label{2.14}
    \left ( \frac {\sum^{n-2}_{i=1}W_i}{\sum^{n-1}_{i=1}W_i} \right )^{\sum^{n-2}_{i=1}W_i}W^{W_{n-1}}_{n-1} \geq \left ( \frac {\sum^{n-1}_{i=1}W_i}{\sum^{n}_{i=1}W_i} \right )^{\sum^{n-1}_{i=1}W_i}W^{W_{n-1}}_{n}.
\end{align}
   The right-hand side expression above when regarded as a function of $W_n$ only is maximized at 
\begin{align*}
  W_n=\frac {W_{n-1}\sum^{n-1}_{i=1}W_i}{\sum^{n-2}_{i=1}W_i}.
\end{align*}
    As the inequality in \eqref{2.14} becomes an equality with this value of $W_n$, we see that inequality \eqref{2.12} does hold for all $n \geq 3$.

   Note that it follows from the arithmetic-geometric mean inequality that
\begin{align*}
   \frac {W_{n-1}}{\sum^{n-1}_{i=1}W_i} \prod^{n-2}_{i=1}\Big(\frac {W_{i+1}}{w_{i+1}}\Big )^{w_{i+1}/W_{n-1}} \leq  \frac {W_{n-1}}{\sum^{n-1}_{i=1}W_i} \left ( \sum^{n-2}_{i=0}\frac {W_{i+1}}{w_{i+1}} \cdot \frac{w_{i+1}}{W_{n-1}} \right )=1.
\end{align*}
   As one checks easily that the above inequality is strict in our case, we see that inequality \eqref{2.13} also holds. We therefore conclude the existence of sequences $\{ w_i \}^n_{i=1}$ satisfying the conditions of Theorem \ref{thm4}.



\end{document}